\RequirePackage{fix-cm}
\documentclass{svjour3}                     
\smartqed  

\usepackage{graphicx}
\usepackage{epstopdf}
\usepackage{amssymb}
\usepackage[cmex10]{amsmath}
\usepackage{algorithm}
\usepackage{algorithmic}
\usepackage{colortbl}
\usepackage[misc]{ifsym}
\usepackage{cite}
\usepackage[justification=centering]{caption}
\usepackage{float}
\usepackage{subfig}
\usepackage{xcolor}

\begin{document}

\title{Event-Triggered Optimal Tracking Control for Strict-Feedback Nonlinear Systems With Non-Affine Nonlinear Faults}
	
	\author{Ling Wang \and Xin Wang \and Ziming Wang}
	\institute{
		L. Wang \and X. Wang\\
		Chongqing Key Laboratory of Nonlinear Circuits and Intelligent Information Processing, College of Electronic and Information Engineering, Southwest University, Chongqing, 400044, P.R. China \\
		Z. M. Wang(\Letter)\\
		The Hong Kong University of Science and Technology (Guangzhou), Nansha, Guangzhou, 511458, P.R. China.\\
		\email{wwwangziming@163.com}\\
	}
	\date{Received: date / Accepted: date}
	\maketitle
	\noindent
	\textbf{Abstract}\\
This article studies the control ideas of the optimal backstepping technique, proposing an event-triggered optimal tracking control scheme for a class of strict-feedback nonlinear systems with non-affine and nonlinear faults. A simplified identifier-critic-actor framework is employed in the reinforcement learning algorithm to achieve optimal control. The identifier estimates the unknown dynamic functions, the critic evaluates the system performance, and the actor implements control actions, enabling modeling and control of anonymous systems for achieving optimal control performance. In this paper, a simplified reinforcement learning algorithm is designed by deriving update rules from the negative gradient of a simple positive function related to the Hamilton-Jacobi-Bellman equation, and it also releases the stringent persistent excitation condition. Then, a fault-tolerant control method is developed by applying filtered signals for controller design. Additionally, to address communication resource reduction, an event-triggered mechanism is employed for designing the actual controller. Finally, the proposed scheme's feasibility is validated through theoretical analysis and simulation.

\keywords{Identifier-critic-actor architecture \and optimized backstepping \and fault-tolerant control \and event-triggered control \and strict-feedback systems}

\section{Introduction\vspace{-0.8em}}
With the problem of resource scarcity day by day, optimal control has received widespread attention\cite{1,2,3,4}. The implementation of optimal control contributes to the reduction of operational costs, enhancement of system efficiency, reinforcement of robustness\cite{5}, and acceleration of response by minimizing unnecessary energy or material waste. It is well known that nonlinear systems are widely present in various fields such as engineering\cite{1111}, physics, and social sciences, so studying the optimal control of nonlinear systems is meaningful. In the context of optimal control problems for nonlinear dynamical systems, solving the Hamilton-Jacobi-Bellman equation (HJBE) is typically employed to derive the optimal control strategy \cite{6}. However, due to the inherent nonlinearity and intractability of the HJBE, obtaining an analytical solution often poses significant challenges, rendering it intractable in numerous cases. Subsequently, Approximate Dynamic Programming (ADP) was introduced\cite{7}. The ADP method has the capability to solve a wider variety of real-world practical issues, such as autonomous driving\cite{8}, robust control\cite{9}, and robot control\cite{10}, greatly expanding the application domains of the ADP method. Reinforcement learning (RL) (or ADP) is a technique applied to train machine learning models to take actions in specific scenarios to maximize expected returns\cite{11}. In recent years, utilizing the consistent approximation and adaptive capability of neural networks (NNs)\cite{2222}, NN-based RL techniques have successfully developed various effective optimal control strategies (such as\cite{13,14,15}). It should be noted that the adaptive backstepping control methods mentioned in references\cite{13} and \cite{14} assume that the system dynamics are known. However, in actual control systems, the dynamics are often unknown, which limits the application of these schemes. To overcome this limitation, observers based on NNs or fuzzy logic systems are commonly used to address the issue of unknown dynamics in nonlinear systems, as indicated in references\cite{3333} and\cite{4444}. However, this traditional control scheme will significantly increase computational complexity. To address this challenge, Wen $et$ $al$.\cite{16} proposed a new optimal backstepping technique. The core idea of this approach is to devise the actual control and all virtual controls as the optimal solutions for the backstepping processes of each subsystem, thereby optimizing the entire backstepping control system. This technique has been widely applied and achieved significant research results (such as\cite{17,18}).

The aforementioned results implicitly indicate the use of a time-triggered control strategy. However, this strategy incurs significant computational and communication resources and struggles to adapt to dynamically changing environments. For applications requiring higher flexibility, efficiency, and responsiveness, event-triggered control(ETC) emerges as a superior choice. ETC technology is a commonly used mode in control systems. It allows for timely response to events occurring in the system, enabling real-time control. Compared to periodic polling methods, ETC only executes relevant operations when an event occurs, avoiding unnecessary resource consumption and improving resource utilization in the system. Therefore, ETC is of significant importance for achieving efficient and reliable control strategies and enhancing the performance of automation systems. In recent years, the field of ETC has witnessed a surge in scholarly interest, with numerous researchers engaging in both theoretical and applied investigations \cite{5555,6666,19,20,21,22,7777,8888}. For instance, reference \cite{5555} focuses on state-constrained inclusive control of nonlinear multi-agent systems using event-triggered inputs. In \cite{6666}, the author designs an event-triggered adaptive inclusive control strategy for heterogeneous stochastic nonlinear multi-agent systems. reference\cite{19} investigates ETC design for heterogeneous multi-agent systems for collaborative control. In\cite{20}, the authors study the influence of bounded disturbances on decentralized ETC systems. Reference\cite{21} discusses how to describe and analyze communication traffic in nonlinear ETC systems using abstract methods. It is worth noting that the results reported in references \cite{5555}-\cite{8888} are based on the premise of fault-free operating conditions.

Many practical systems, particularly those in industrial domains, exhibit the characteristics of non-affine and nonlinear faults. These fault systems often encounter various challenges, such as sensor failures and actuator faults. Therefore, it is of paramount importance to develop a fault-tolerant control (FTC) scheme that possesses both effective fault tolerance capabilities and enhances system availability and cost-effectiveness. Recently, there have been numerous outstanding FTC solutions developed in the literature\cite{23,24,25,26,27,28}. For example, reference\cite{23} addresses the problem of FTC for attitude stabilization under multiple disturbances. Shen $et$ $al$. achieve integral sliding mode FTC for spacecraft attitude stability in reference\cite{24}. Although significant progress has been made in dealing with actuator faults in the aforementioned research, to the best of the authors' knowledge, there is limited literature on the event-triggered optimal FTC problem for strict-feedback nonlinear systems(SFNSs) with unknown dynamic functions. Addressing this issue requires tackling the following three challenging problems: 1) How to develop an optimal tracking control scheme for SFNSs under unknown dynamics, non-affine nonlinear faults, and acceptable complexity? 2) How can resource-saving strategies be introduced into control systems to effectively improve resource utilization efficiency? 3) How can we ensure that the controller achieves satisfactory tracking performance while operating under the FTC strategy and event-triggered mechanism? These challenges motivate our research in this area.

In conclusion, a fresh optimal FTC scheme is put forward based on an identifier-critic-actor architecture and RL for unknown dynamic SFNSs. The summary article makes three main contributions, which are as follows: 1) The offered optimal control is remarkably simplified, as a simple positive function associated with the HJBE is utilized to calculate the negative gradient, which in turn determines the updating rate of RL. Moreover, constructing an identifier-critic-actor architecture provides a comprehensive and effective control solution. Therefore, the simplified structure and extensive applicability of this algorithm make it easier to execute and promote. Note that the past approaches mainly trained RL models by gradient descent on the square of HJBE, however, this algorithm is very complex and tricky and requires persistent excitation. These defects limit the application and scalability of the approach. Therefore, simplifying the optimal control approach becomes more meaningful. 2) In modern control systems, non-affine nonlinear faults are widely present, which may lead to a decrease in system performance or even system failure. The use of an approach based on Butterworth low-pass filters and neural network approximation in practical controller design to develop adaptive tracking FTC laws for systems can effectively enhance system fault tolerance and stability, thus compensating for the adverse effects of non-affine nonlinear faults. 3) By incorporating an event-triggered mechanism, the update rate of the practical controller is further adjusted to improve the efficiency of communication resource utilization. Designing an event-triggered optimal controller for a SFNS under unknown dynamics and non-affine nonlinear faults is a challenge that has not been addressed in published works. It is worth noting that although Pang et al. \cite{5555} considered ETC, the system omitted non-affine nonlinear faults.

\section{Preliminaries\vspace{-0.8em}}
\subsection{Problem Formulation \vspace{-0.8em}}
Consider the following system:
\begin{align}
\dot{x}_i & =x_{i+1}+f_i\left(\bar{x}_i\right) \quad {i}=1, \ldots, {n}-1 \notag\\
\dot{x}_n & =u+f_n\left(\bar{x}_n\right)+\sigma\left(t-T_0\right) \lambda(x, u) \quad {n} \geq 2\notag \\
{y} & =x_1
\end{align}
where $\bar{x}_i=\left[x_1, \ldots, x_i\right]^T  \in \mathbb{R}^i$ is the state vector and $x_1 \in \mathbb{R}$ is the system output. Besides, ${u} \in \mathbb{R}$ and  $f_i\left(\bar{x}_i\right) \in \mathbb{R}$ represent the control input and the unknown smooth nonlinear function, respectively. $\lambda(x, u)$ is an unknown disturbance caused by a fault, and $\sigma\left(t-T_0\right)$ represents a time curve of a fault occurring at an unspecified time $T_0$, expressed as
$$
\sigma\left(t-T_0\right)= \begin{cases}0, & t<T_0 \\ 1-\mathrm{e}^{-\alpha\left(t-T_0\right)}, & t \geq T_0\end{cases}
$$
where $\alpha>0$ signifies the rate at which the unknown fault progresses.

$Remark\hspace{0.5em}\textsl{1:}$\label{Remark1}
In fact, nonlinear systems are prevalent in practical engineering applications, and the majority of them exhibit non-affine nonlinear faults, such as electromechanical systems, continuous flow pendulum systems, and one-link manipulator system. They can all be modeled as system (1).

Here is a definition of admissible control. For more detailed information on the field of optimal control, please refer to references \cite{16,17}.

$Definition\hspace{0.5em}\textsl{1\cite{16} :}$\label{Definition1}
Let $\dot{{x}}({t})={f}({x})+{g}({x}) {u}({x})$ be a nonlinear system, the control protocol ${u}({x})$ of the nonlinear system is said to be admissible on the set $\mathcal{W}$, which is denoted by ${u}({x}) \in \Psi(\mathcal{W})$, if ${u}({x})$ is continous with ${u}(0)=0$, and stabilizes the nonlinear system on $\mathcal{W}$, and makes the infinite horizon value function ${J}({x})=\int_t^{\infty} c(x(s), u(x)) d s$ finite.

Optimal Control: The optimal control objective of the system $\dot{{x}}({t})={f}({x})+{g}({x}) {u}({x})$ is to find an admissible control that minimizes the performance index function ${J}({x})$ in order to achieve optimal system performance.

$Lemma\hspace{0.5em}\textsl{1\cite{17} :}$\label{lemma1}
If, for a given function ${V}$, its derivative $\dot{V}$ meets $\dot{V}\leq-a V+m$, where $a>0$ and $m>0$, thus, the ensuing inequality is valid:
\begin{equation}
{V} \leq e^{-a t} V(0)+\frac{m}{a}\left(1-e^{-a t}\right).
\end{equation}

\newtheorem{assumption}{Assumption}
$Assumption\hspace{0.5em}\textsl{1\cite{18} :}$\label{assumption1}	
An unknown non-negative function $\bar{h}(x, u)$ lives, fulfilling the following inequality:
\begin{equation}
\left|f_n\left(\bar{x}_n\right)+\sigma\left(t-T_0\right) \lambda(x, u)\right| \leq \bar{h}(x, u).
\end{equation}

Control Objective: This paper proposes an optimal FTC strategy based on the optimal backstepping control technique for a kind of SFNSs with non-affine nonlinear faults. Additionally, an event-triggered procedure is engineered to further minimize the utilization of resources. As a result, all error signals are ensured to be semiglobally uniformly ultimately bounded (SGUUB), possess excellent tracking performance, and avoid Zeno behavior.

\section{Main Results\vspace{-0.8em}}
The focus of this section lies in the introduction of an adaptive optimal tracking control algorithm that utilizes the identifier-actor-critic structure.
\subsection{Event-Triggered Optimal Backstepping Control\vspace{-0.8em}}
Before proceeding, given below are the state coordinate transformations:
\begin{align}
e_1& =x_1-y_r \\
{e_i}& =x_i-\hat{\alpha}_{i-1}^*,\quad {i}=2, \ldots, {n} 
\end{align}
where $\hat{\alpha}_{i-1}^*$ represents the optimal virtual control will be designed in step $i - 1$ and $y_r$ denotes the provided reference signal.

\emph{Step \textsl{1 :}} Combining (1) and (4), the time derivative of ${e}_1$ is
\begin{equation}
\dot{e}_1=x_2+f_1\left(\bar{x}_1\right)-\dot{y}_r.
\end{equation}

Construct the performance index function as
\begin{equation}
\begin{aligned}
J_1^*(e_1) & =\min _{\alpha_1 \in \Psi(\mathcal{W})}\left(\int_t^{\infty} c_1\left(e_1(s), \alpha_1\left(e_1\right)\right) d s\right) \\
& =\int_t^{\infty} c_1\left(e_1(s), \alpha_1^*\left(e_1\right)\right) d s,
\end{aligned}
\end{equation}
where $c_1\left(e_1, \alpha_1\right)=e_1^2+\alpha_1^2$ denotes the cost function, $\alpha_1^*$ and $\mathcal{W}$ represent the optimal virtual control and a compact set containing origin, respectively.

Considering $x_2 \triangleq \alpha_1^*\left(e_1\right)$, we can obtain the HJBE by calculating the temporal derivatives of the optimal performance function, resulting in the following form:
\begin{equation}
\begin{aligned}
H_1\left(e_1, \alpha_1^*, \frac{d J_1^*\left(e_1\right)}{d e_1}\right)
=& e_1^2+\alpha_1^{* 2}+\frac{d J_1^*\left(e_1\right)}{d e_1} \times\left(x_2+f_1\left(\bar{x}_1\right)-\dot{y}_r\right)\\
=&0.
\end{aligned}
\end{equation}

By resolving $\frac{\partial \mathrm{H_1}}{\partial \alpha_1^*}=0$, we have
\begin{equation}
\alpha_1^*=-\frac{1}{2} \frac{d J_1^*\left(e_1\right)}{d e_1}.
\end{equation}

To achieve the expected outcome, the term $\frac{d J_1^*\left(e_1\right)}{d e_1}$ is decomposed as
\begin{equation}
\frac{d J_1^*\left(e_1\right)}{d e_1}=2 \rho_1 e_1+2 f_1\left(\bar{x}_1\right)+J_1^0\left(\bar{x}_1, e_1\right),
\end{equation}
where $\rho_1>0$ is a design constant, and $J_1^0\left(\bar{x}_1, e_1\right) \in \mathbb{R}$ is a continuous function defined as $
J_1^0\left(\bar{x}_1, e_1\right)=-2 \rho_1 e_1-2 f_1\left(\bar{x}_1\right)+\frac{d J_1^*\left(e_1\right)}{d e_1}$.

The substitution of (10) into (9) yields the subsequent expression
\begin{equation}
\alpha_1^*=-\rho_1 e_1-f_1\left(\bar{x}_1\right)-\frac{1}{2} J_1^0\left(\bar{x}_1, e_1\right).
\end{equation}

Given the capability of NNs to approximate any continuous function, we can approximate the unknown continuous functions $f_1\left(\bar{x}_1\right)$ and $J_1^0\left(\bar{x}_1, e_1\right)$ as
\begin{align}
f_1\left(\bar{x}_1\right) & =\S_{f_1}^{* T} E_{f_1}\left(\bar{x}_1\right)+\omega_{f_1}\left(\bar{x}_1\right) \\
J_1^0\left(\bar{x}_1, e_1\right) & =\S_{J_1}^{* T} E_{J_1}\left(\bar{x}_1, e_1\right)+\omega_{J_1}\left(\bar{x}_1, e_1\right)
\end{align}
where $\S_{f_1}^*$ and $\S_{J_1}^*$ represent the ideal weights, $E_{f_1}$ and $E_{J_1}$ are denoted as the basis function vectors and the NN approximation errors are represented by $\omega_{f_1}$ and $\omega_{J_1}$.

Substituting (12) and (13) into (10) and (11), one has
\begin{align}
\alpha_1^*=&-\rho_1 e_1-\frac{1}{2} \omega_1-\frac{1}{2} \S_{J_1}^{* T} E_{J_1}\left(\bar{x}_1, e_1\right)-\S_{f_1}^{* T} E_{f_1}\left(\bar{x}_1\right)  \\
\frac{d J_1^*\left(e_1\right)}{d e_1}= &2 \rho_1 e_1+\omega_1+\S_{J_1}^{* T} E_{J_1}\left(\bar{x}_1, e_1\right)+2 \S_{f_1}^{* T} E_{f_1}\left(\bar{x}_1\right) 
\end{align}
where $\omega_1=2 \omega_{f_1}+\omega_{J_1}$.

Utilizing the optimal virtual controller directly is not feasible due to the unknown ideal weight vectors $\S_{f_1}^*$ and $\S_{J_1}^*$, Therefore, an identifier-critic-actor structure is employed based on NN approximation algorithms to achieve the desired results.

The following represents the NNs for the identifier, critic, and actor used for approximating unknown dynamic functions, assessing control performance, and effectuating control actions:
\begin{align}
\hat{f}_1\left(\bar{x}_1\right)=& \hat{\S}_{f_1}^T E_{f_1}\left(\bar{x}_1\right) \\
\frac{\mathrm{d} \hat{J}_1^*\left(e_1\right)}{\mathrm{d} e_1}=& 2 \rho_1 e_1+2 \hat{\S}_{f_1}^T E_{f_1}\left(\bar{x}_1\right)+\hat{\S}_{c_1}^T E_{J_1}\left(\bar{x}_1, e_1\right) \\
\hat{\alpha}_1^*=& -\rho_1 e_1-\hat{\S}_{f_1}^T E_{f_1}\left(\bar{x}_1\right)-\frac{1}{2} \hat{\S}_{a_1}^T E_{J_1}\left(\bar{x}_1, e_1\right)
\end{align}
where $\hat{f}_1\left(\bar{x}_1\right)$, $\frac{\mathrm{d} \hat{J}_1^*\left(e_1\right)}{\mathrm{d} e_1}$, and $\hat{\alpha}_1^*$ are the identifier output, the estimation of $\frac{d J_1^*\left(e_1\right)}{d e_1}$ and $\alpha_1^*$, respectively. The weights of the identifier, critic, and actor NNs are represented by $\hat{\S}_{f_1}$, $\hat{\S}_{c_1}$, and $\hat{\S}_{a_1}$,
respectively. In the identifier-critic-actor model, the updating rules for weight vectors $\dot{\hat{\S}}_{f_1}, \dot{\hat{\S}}_{c_1}, \dot{\hat{\S}}_{a_1}$ are designed as
\begin{align}
\dot{\hat{\S}}_{f_1}= & \Pi_1\left(E_{f_1}\left(\bar{x}_1\right) e_1-\gamma_1 \hat{\S}_{f_1}\right) \\
\dot{\hat{\S}}_{c_1}= & -\varepsilon_{c_1} E_{J_1}\left(\bar{x}_1, e_1\right) E^T_{J_1}\left(\bar{x}_1, e_1\right) \hat{\S}_{c_1} \\
\dot{\hat{\S}}_{a_1}= & -E_{J_1}\left(\bar{x}_1, e_1\right) E^T_{J_1}{ }\left(\bar{x}_1, e_1\right)\times\left(\varepsilon_{a_1}\left(\hat{\S}_{a_1}-\hat{\S}_{c_1}\right)+\varepsilon_{c_1} \hat{\S}_{c_1}\right)
\end{align}
where $\Pi_1$ is a positive-definite matrix and $\gamma_1$, $\varepsilon_{c_1}$, $\varepsilon_{a_1}$, $\rho_1$ are all design parameters, which are selected to satisfy $\gamma_1>0$, $\varepsilon_{a_1}>\frac{1}{2}$, $\varepsilon_{a_1}>\varepsilon_{c_1}>\frac{\varepsilon_{a_1}}{2}$, $\rho_1>3$.

\emph{Remark \textsl{2:}} The RL update rules outlined above are derived from the negative gradient of a straightforward, positively correlated function associated with the HJBE. The following details will be provided on how to derive the relevant parameters $\dot{\hat{\S}}_{c_1}$ and $\dot{\hat{\S}}_{a_1}$.

Substituting (17) and (18) into (8), the HJBE is obtained as
\begin{equation}
\begin{aligned}
H_1\left(e_1, \hat{\alpha}_1^*, \frac{d \hat{J}_1^*\left(e_1\right)}{d e_1}\right)=&e_1^2+\left(\rho_1 e_1+\hat{\S}_{f_1}^T E_{f_1}\left(\bar{x}_1\right)+\frac{1}{2} \hat{\S}_{a_1}^T E_{J_1}\left(\bar{x}_1, e_1\right)\right)^2 \\
&+\left(2 \rho_1 e_1+2 \hat{\S}_{f_1}^T E_{f_1}\left(\bar{x}_1\right)+\hat{\S}_{c_1}^T E_{J_1}\left(\bar{x}_1, e_1\right)\right) \\
&\times\left(-\rho_1 e_1+f_1\left(\bar{x}_1\right)-\hat{\S}_{f_1}^T E_{f_1}\left(\bar{x}_1\right)-\frac{1}{2} \hat{\S}_{a_1}^T E_{J_1}\left(\bar{x}_1, e_1\right)-\dot{y}_r\right).
\end{aligned}
\end{equation}

The Bellman residual $\Xi_1$ is defined as follows
\begin{equation}
\begin{aligned}
\Xi_1& =H_1\left(e_1, \hat{\alpha}_1^*, \frac{\mathrm{d} \hat{J}_1^*\left(e_1\right)}{\mathrm{d} e_1}\right)-H_1\left(e_1, \alpha_1^*, \frac{d J_1^*(e_1)}{d e_1}\right) \\
& =H_1\left(e_1, \hat{\alpha}_1^*, \frac{\mathrm{d} \hat{J}_1^*\left(e_1\right)}{\mathrm{d} e_1}\right).
\end{aligned}
\end{equation}

Relying on the previous analysis, it can be concluded that the expected value of the optimized solution $\hat{\alpha}_1^*$ satisfies $\Xi_1 \rightarrow 0$. If $H_1\left(e_1, \hat{\alpha}_1^*, \frac{\mathrm{d} \hat{J}_1^*\left(e_1\right)}{\mathrm{d} e_1}\right)=0$ holds and there exists a unique solution, it can be equivalently expressed as the following equation:
\begin{equation}
\frac{\partial H_1\left(e_1, \hat{\alpha}_1^*, \frac{d \hat{J}_1^*\left(e_1\right)}{d e_1}\right)}{\partial \hat{\S}_{a_1}}=\frac{1}{2} E_{J_1} E_{J_1}^T\left(\hat{\S}_{a_1}-\hat{\S}_{c_1}\right)=0.
\end{equation}

To ensure that the derived RL update rate satisfies (24), we construct the positive function as follows: $K_1=\left(\hat{\S}_{a_1}-\hat{\S}_{c_1}\right)^T\left(\hat{\S}_{a_1}-\hat{\S}_{c_1}\right)$. It is evident from (24) that $K_1=0$ holds. From $\frac{\partial K_1}{\partial \hat{\S}_{a_1}}=-\frac{\partial K_1}{\partial \hat{\S}_{c_1}}=2\left(\hat{\S}_{a_1}-\hat{\S}_{c_1}\right)$, we can obtain
\begin{equation}
\begin{aligned}
\frac{{d} K_1}{{dt}}& =\frac{\partial K_1}{\partial \hat{\S}_{c_1}^T} \times \dot{\hat{\S}}_{c_1}+\frac{\partial K_1}{\partial \hat{\S}_{a_1}^T} \times \dot{\hat{\S}}_{a_1} \\
& =-\frac{\varepsilon_{a_1}}{2} \frac{\partial K_1}{\partial \hat{\S}_{a_1}^T} E_{J_1}\left(\bar{x}_1, e_1\right) E_{J_1}^T\left(\bar{x}_1, e_1\right) \frac{\partial K_1}{\partial \hat{\S}_{a_1}} \\
& \leq 0.
\end{aligned}
\end{equation}

Equation (25) indicates that using the update rates (20) and (21) ensures the validity of $K_1 \rightarrow 0$.

\emph{Remark \textsl{3:}} To illustrate the algorithm's simplification, a comparison with the approach outlined in \cite{16} is presented below. 

The critic and actor in \cite{16} both employ updated rules of the form presented below
\begin{align}
\dot{\hat{\S}}_{c_1} =&-\frac{\varepsilon_{c 1}}{\left\|\omega_1\right\|^2+1} \omega_1\left(\omega_1^T \hat{\S}_{c_1}-\left(\rho_1^2-1\right)\right.{e}_1^2  +2 \rho_1 {e}_1\left(f_1\left(x_1\right)-\dot{y}_r\right)
+\frac{1}{4} \hat{\S}_{a_1}^T \left. E_{J_1} E_{J_1}^T \hat{\S}_{a_1}\right)\notag\\
\dot{\hat{\S}}_{a_1}=&\frac{1}{2} E_{J_1} {e}_1-\varepsilon_{a_1} E_{J_1} E_{J_1}^T \hat{\S}_{a_1} +\frac{\varepsilon_{c_1}}{4\left(\left\|\omega_1\right\|^2+1\right)} E_{J_1}E_{J_1}^T \hat{\S}_{a_1} \omega_1^T \hat{\S}_{c_1}\notag
\end{align}
where $\varepsilon_{c_1}>0, \varepsilon_{a_1}>0, \omega_1=E_{J_1}\left(f_1\left(x_1\right)-\rho_1 {e}_1-\right.$ $\left.\frac{1}{2} \hat{\S}_{a_1}^T E_{J_1}-\dot{y}_r\right)$. Upon comparing with (20) and (21), it becomes evident that the proposed optimized control is algorithmically simpler. Furthermore, achieving $\omega_1\neq0$ also requires effective training of adaptive parameters. Therefore, the requirement for persistence excitation, i.e.,  $\eta_1 I_{m_1} \leq \omega_i \omega_i^T \leq \zeta_1 I_{m_1}$ where $\eta_1\textgreater0$ and $\zeta_1\textgreater0$ are constants, is essential. Therefore, it can be concluded that the optimized control proposed in this paper releases the condition of persistent excitation.

\emph{Step \textsl{i ( i = 2,...,n-1) :}}
Similarly, according to (1) and (5), we can obtain
\begin{equation}
\dot{e}_i=x_{i+1}+f_i\left(\bar{x}_i\right)-\dot{\hat{\alpha}}_{i-1}^*.
\end{equation}

The performance function can be described as
\begin{equation}
\begin{aligned}
J_i^*(e_i) & =\min _{\alpha_i \in \Psi(\mathcal{W})}\left(\int_t^{\infty} c_i\left(e_i(s), \alpha_i\left(e_i\right)\right) d s\right) \\
& =\int_t^{\infty} c_i\left(e_i(s), \alpha_i^*\left(e_i\right)\right) d s,
\end{aligned}
\end{equation}
where $c\left(e_i, \alpha_i\right)=e_i^2+\alpha_i^2$ denotes the cost function, and $\alpha_i^*$ represents the optimal virtual control.

Considering $x_{i+1} \triangleq \alpha_i^*\left(e_i\right)$, the HJBE can be yielded as follows
\begin{equation}
\begin{aligned}
H_i\left(e_i, \alpha_i^*, \frac{d J_i^*\left(e_i\right)}{d e_i}\right) =& e_i^2+\alpha_i^{* 2}+\frac{d J_i^*\left(e_i\right)}{d e_i} \times\left(\alpha_i^*\left(e_i\right)+f_i\left(\bar{x}_i\right)- \dot{\hat{\alpha}}_{i-1}^*\right) \\
=& 0.
\end{aligned}
\end{equation}

Then $\alpha_i^*$ can be acquired through resolving $\frac{\partial H_i}{\partial \alpha_i^*}=0$ as
\begin{equation}
\alpha_i^*=-\frac{1}{2} \frac{d J_i^*\left(e_i\right)}{d e_i}.
\end{equation}

Similar to step 1, the term $\frac{d J_i^*\left(e_i\right)}{d e_i}$ is decomposed as
\begin{equation}
\frac{d J_i^*\left(e_i\right)}{d e_i}=2 \rho_i e_i+2 f_i\left(\bar{x}_i\right)+J_i^0\left(\bar{x}_i, e_i\right),
\end{equation}
where $\rho_i>0$, and $J_i^0\left(\bar{x}_i, e_i\right)$ is a continuous function, which can be yielded as
$J_i^0\left(\bar{x}_i, e_i\right)=-2 \rho_i e_i-2 f_i\left(\bar{x}_i\right)+\frac{d J_i^*\left(e_i\right)}{d e_i}
$.

By substituting (30) into (29), we have
\begin{equation}
\alpha_i^*=-\rho_i e_i-f_i\left(\bar{x}_i\right)-\frac{1}{2} J_i^0\left(\bar{x}_i, e_i\right).
\end{equation}

Similarly, we can obtain the following NNs for the identifier, critic, and actor:
\begin{align}
\hat{f}_i\left(\bar{x}_i\right)=&\hat{\S}_{f_i}^T E_{f_i}\left(\bar{x}_i\right) \\
\frac{\mathrm{d} \hat{J}_i^*\left(e_i\right)}{\mathrm{d} e_i}=&2 \rho_i e_i+2 \hat{\S}_{f_i}^T E_{f_i}\left(\bar{x}_i\right)+\hat{\S}_{c_i}^T E_{J_i}\left(\bar{x}_i, e_i\right)\\
\hat{\alpha}_i^*=&-\rho_i e_i-\hat{\S}_{f_i}^T E_{f_i}\left(\bar{x}_i\right)-\frac{1}{2} \hat{\S}_{a_i}^T E_{J_i}\left(\bar{x}_i, e_i\right)
\end{align}
where $\hat{f}_i\left(\bar{x}_i\right)$, $\frac{\mathrm{d} \hat{J}_i^*\left(e_i\right)}{\mathrm{d} e_i}$, and $\hat{\alpha}_i^*$ are the identifier output, the estimation of $\frac{d J_i^*\left(e_i\right)}{d e_i}$ and $\alpha_i^*$,
respectively. The weights of the identifier, critic, and actor NNs are represented by $\hat{\S}_{f_i}$, $\hat{\S}_{c_i}$, and $\hat{\S}_{a_i}$,
respectively. The rules for updating these weights are formulated as
\begin{align}
\dot{\hat{\S}}_{f_i}=&\Pi_i\left(E_{f_i}\left(\bar{x}_i\right) e_i-\gamma_i \hat{\S}_{f_i}\right) \\
\dot{\hat{\S}}_{c_i}=&-\varepsilon_{c_i} E_{J_i}\left(\bar{x}_i, e_i\right) E^T_{J_i}{ }\left(\bar{x}_i, e_i\right) \hat{\S}_{c_i} \\
\dot{\hat{\S}}_{a_i}=&-E_{J_i}\left(\bar{x}_i, e_i\right) E^T_{J_i}{ }\left(\bar{x}_i, e_i\right) \times \left(\varepsilon_{a_i}\left(\hat{\S}_{a_i}-\hat{\S}_{c_i}\right)+\varepsilon_{c_i} \hat{\S}_{c_i}\right)
\end{align}
where $\Pi_i$ is a positive-definite matrix and $\gamma_i$, $\varepsilon_{c_i}$, $\varepsilon_{a_i}$, $\rho_i$ are all design parameters, which are chosen to satisfy $\gamma_i>0$, $\varepsilon_{a_i}>\frac{1}{2}$, $\varepsilon_{a_i}>\varepsilon_{c_i}>\frac{\varepsilon_{a_i}}{2}$, $\rho_i>3$.

\emph{Step \textsl{n :}}
Combining equation (1) and (5), we have
\begin{equation}
\dot{e}_n=u+f_n\left(\bar{x}_n, u\right)+\sigma\left(t-T_0\right) \lambda(x, u)-\dot{\hat{\alpha}}_{n-1}^*.
\end{equation}

Let $u^*$ denote the optimal control. Then the performance function can be expressed as
\begin{equation}
\begin{aligned}
J_n^*\left(e_n\right)& =\min _{u \in \Psi(\mathcal{W})}\left(\int_t^{\infty} c_n\left(e_n(s), u\left(e_n\right)\right) d s\right) \\
& =\int_t^{\infty} c_n\left(e_n(s), u^*\left(e_n\right)\right) d s,
\end{aligned}
\end{equation}
where $c_n\left(e_n, u\right)=e_n^2+u^2$ denotes the cost function. 

Combining equation (38), the HJBE is defined as
\begin{align}
H_n\left(e_n, u^*, \frac{d J_n^*\left(e_n\right)}{d e_n}\right)=& e_n^2+u^{* 2}+\frac{d J_n^*\left(e_n\right)}{d e_n} \times\bigg(u^*+f_n\left(\bar{x}_n, u\right)-\dot{\hat{\alpha}}_{n-1}^* +\sigma(t-T_0) \lambda(x, u)\bigg)\notag\\
=& 0.
\end{align}

From $\frac{\partial H_n}{\partial u^*}=0$, one has
\begin{equation}
u^*=-\frac{1}{2}\frac{d J_n^*\left(e_n\right)}{d e_n}.
\end{equation}

The term $\frac{d J_n^*\left(e_n\right)}{d e_n}$ is decomposed as
\begin{equation}
\frac{d J_n^*\left(e_n\right)}{d e_n}=2 \rho_n e_n+2 f_n\left(\bar{x}_n, u\right)+J_n^0\left(\bar{x}_n, e_n\right),
\end{equation}
where $\rho_n>0$, and $J_n^0\left(\bar{x}_n, e_n\right)$ is a continuous function, which can be yielded as
$J_n^0\left(\bar{x}_n, e_n\right)=-2 \rho_n e_n-2 f_n\left(\bar{x}_n,u\right)+\frac{d J_n^*\left(e_n\right)}{d e_n}
$.

Substituting (42) into (41), one obtains
\begin{equation}
u^*=-\rho_n e_n-f_n\left(\bar{x}_n, u\right)-\frac{1}{2} J_n^0\left(\bar{x}_n, e_n\right).
\end{equation}

Smiliar to (12) and (13), we can approximate the unknown continuous functions $f_n\left(\bar{x}_n\right)$ and $J_n^0\left(\bar{x}_n, e_n\right)$ as
\begin{align}
f_n\left(\bar{x}_n, u\right)&=\S_{f_n}^{* T} E_{f_n}\left(\bar{x}_n, u\right)+\omega_{f_n}\left(\bar{x}_n\right)\\
J_n^0\left(\bar{x}_n, e_n\right)&=\S_{J_n}^{* T} E_{J_n}\left(\bar{x}_n, e_n\right)+\omega_{J_n}\left(\bar{x}_n, e_n\right)
\end{align}
where $\S_{f_n}^*$ and $\S_{J_n}^*$ are the ideal weights, $E_{f_n} $ and $E_{J_n}$ are denoted as the basis function vectors, and the NN approximation errors are represented by $\omega_{f_n}$ and $\omega_{J_n}$.

In the same manner, as before, we can derive the following NNs for the identifier, critic, and actor:
\begin{align}
\hat{f}_n\left(\bar{x}_n, u\right)=&\hat{\S}_{f_n}^T E_{f_n}\left(\bar{x}_n, u\right) \\
\frac{d\hat{J}_n^*\left(e_n\right)}{d e_n}=&2 \rho_n e_n+2 \hat{\S}_{f_n}^T E_{f_n}\left(\bar{x}_n, u\right)+\hat{\S}_{c_n}^T E_{J_n}\left(\bar{x}_n, e_n\right)\\
\hat{u}=&-\rho_n e_n-\hat{\S}_{f_n}^T E_{f_n}\left(\bar{x}_n, u\right)-\frac{1}{2} \hat{\S}_{a_n}^T E_{J_n}\left(\bar{x}_n, e_n\right)
\end{align}
where $\hat{f}_n\left(\bar{x}_n, u\right), \frac{d\hat{J}_n^*\left(e_n\right)}{d e_n}$, and $\hat{u}^*$ are the identifier output, the estimation of $\frac{d J_n^*\left(e_n\right)}{d e_n}$ and $u^*$, respectively. The weights of the identifier, critic, and actor NNs are represented by $\hat{\S}_{f_n}$, $\hat{\S}_{c_n}$, and $\hat{\S}_{a_n}$,
respectively. The rules for updating these weights are devised as
\begin{align}
\dot{\hat{\S}}_{f_n}=&\Pi_n\left(E_{f_n}\left(\bar{x}_n, u\right) e_n-\gamma_n \hat{\S}_{f_n}\right)\\
\dot{\hat{\S}}_{c_n}=&-\varepsilon_{c_n} E_{J_n}\left(\bar{x}_n, e_n\right) E^T_{J_n}\left(\bar{x}_n, e_n\right) \hat{\S}_{c_n}\\
\dot{\hat{\S}}_{a_n}=&-E_{J_n}\left(\bar{x}_n, e_n\right) E^T_{J_n}{\left(\bar{x}_n, e_n\right)}\times\left(\varepsilon_{a_n}\left(\hat{\S}_{a_n}-\hat{\S}_{c_n}\right)+\varepsilon_{c_n} \hat{\S}_{c_n}\right)
\end{align}
where $\Pi_n$ is a positive-definite matrix and $\gamma_n$, $\varepsilon_{c_n}$, $\varepsilon_{a_n}$, $\rho_n$ are all design parameters, which satisfy $\gamma_n>0$, $\varepsilon_{a_n}>\frac{1}{2}$, $\varepsilon_{a_n}>\varepsilon_{c_n}>\frac{\varepsilon_{a_n}}{2}$, $\rho_n>3$.

To conserve communication resources, combining equation (48), the controller $u({t})$ is designed based on an event-triggered approach as the following forms:
\begin{align}
U({t})& =-(1+\beta)\left(\hat{u} \tanh \left(\frac{e_n \hat{u}}{v}\right)+\zeta \tanh \left(\frac{e_n \zeta}{v}\right)\right) \notag\\
u({t})& ={U}\left({t}_k\right), \forall {t} \in\left[{t}_k, {t}_{k+1}\right)
\end{align}
this results in the trigger condition being expressed as
\begin{equation}
\begin{aligned}
|{u}({t})-{U}({t})|<\beta|{u}({t})|+\theta \\
\end{aligned}
\end{equation}
\begin{equation}
\begin{aligned}
t=\left\{\begin{array}{lr}
t_{k+1}, & |u(t)-{U}(t)| \geq \beta|u(t)|+\theta \\
t_k, & \text { else }
\end{array}\right.
\end{aligned}
\end{equation}
where $0<\beta<1, \quad \theta>0, \quad \zeta>\frac{\theta}{1-\beta} \quad$ and $\quad v>0 \quad$ are designed constants. The update of the event trigger time ${t}_k\left({k} \in Z^{+}\right)$ is contingent upon the fulfillment of the trigger condition. When inequality (53) is satisfied, the actual control input ${u}={U}\left({t}_{k+1}\right)$ is updated.

\subsection{Stability Analysis\vspace{-0.8em}}
A theoretical conclusion is derived in this section by applying the relevant conditions of 3.1 and employing the Lyapunov stability analysis method.

We define the corresponding error as $\tilde{\S}_{z_i}=\hat{\S}_{z_i}-\S_{z_i}^*(z=f, c, a)$ for ${i}=1, \ldots, {n}$.

\emph{Step \textsl{1 :}} We choose the following Lyapunov candidate function:
\begin{equation}
\begin{aligned}
V_1=&\frac{1}{2} e_1^2+\frac{1}{2} \tilde{\S}_{f_1}^T \Pi_1^{-1} \tilde{\S}_{f_1}+\frac{1}{2} \tilde{\S}_{c_1}^T \tilde{\S}_{c_1}+\frac{1}{2} \tilde{\S}_{a_1}^T \tilde{\S}_{a_1}.
\end{aligned}
\end{equation}

Letting $e_2=x_2-\hat{\alpha}_1^*$, then (6) can be transformed into
\begin{equation}
\dot{e}_1=\hat{\alpha}_1^*+e_2+f_1\left(\bar{x}_1\right)-\dot{y}_r.
\end{equation}

From equations (12), (18)-(21), and (56), we get
\begin{align}
\begin{aligned}
\dot{V}_1=&e_1\bigg(-\rho_1 e_1-\tilde{\S}_{f_1}^T E_{f_1}-\frac{1}{2} \hat{\S}_{a_1}^T E_{J_1}+e_2 -\dot{y}_r+\omega_{f_1}\bigg)+\tilde{\S}_{f_1}^T\left(E_{f_1} e_1-\gamma_1 \hat{\S}_{f_1}\right) \\
& -\varepsilon_{c_1} \tilde{\S}_{c_1}^T E_{J_1} E_{J_1}^T \hat{\S}_{c_1}-\tilde{\S}_{a_1}^T E_{J_1} E_{J_1}^T \times\left(\varepsilon_{a_1}\left(\hat{\S}_{a_1}-\hat{\S}_{c_1}\right)+\varepsilon_{c_1} \hat{\S}_{c_1}\right).
\end{aligned}
\end{align}

By applying Young's inequality, it yields
\begin{align}
e_1 \omega_{f_1} \leq& \frac{1}{2} e_1^2+\frac{1}{2} \omega_{f_1}^2 \notag\\
e_1 e_2 \leq& \frac{1}{2} e_1^2+\frac{1}{2} e_2^2\notag \\
-e_1\dot{y}_r \leq& \frac{1}{2} e_1^2+\frac{1}{2} \dot{y}_r^2\notag \\
-\frac{1}{2} e_1 \hat{\S}_{a_1}^T E_{J_1}  \leq& \frac{1}{4}\left(\hat{\S}_{a_1}^T E_{J_1}\right)^2 +\frac{1}{4} e_1^2\\
\left(\varepsilon_{a_1}-\varepsilon_{c_1}\right) \tilde{\S}_{a_1}^T E_{J_1} E_{J_1}^T \hat{\S}_{c_1}  \leq& \frac{\left(\varepsilon_{a_1}-\varepsilon_{c_1}\right)}{2}\left(\tilde{\S}_{a_1}^T E_{J_1}\right)^2 +\frac{\left(\varepsilon_{a_1}-\varepsilon_{c_1}\right)}{2}\left(\hat{\S}_{c_1}^T E_{J_1}\right)^2
\end{align}

Based on the fact that $\tilde{\S}_{z_i}=\hat{\S}_{z_i}-\S_{z_i}^*(z=f, c, a)$, one can have
\begin{align}
\tilde{\S}_{f_1}^T \hat{\S}_{f_1}=&\frac{1}{2} \tilde{\S}_{f_1}^T \tilde{\S}_{f_1}+\frac{1}{2} \hat{\S}_{f_1}^T \hat{\S}_{f_1}-\frac{1}{2} \S_{f_1}^{* T} \S_{f_1}^*\\
\tilde{\S}_{c_1}^T E_{J_1} E_{J_1}^T\hat{\S}_{c_1}=& \frac{1}{2}\left(\tilde{\S}_{c_1}^T E_{J_1}\right)^2+\frac{1}{2}\left(\hat{\S}_{c_1}^T E_{J_1}\right)^2  -\frac{1}{2}\left(\S_{J_1}^{* T} E_{J_1}\right)^2 \\
\tilde{\S}_{a_1}^T E_{J_1} E_{J_1}^T\hat{\S}_{a_1} =&\frac{1}{2}\left(\tilde{\S}_{a_1}^T E_{J_1}\right)^2+\frac{1}{2}\left(\hat{\S}_{a_1}^T E_{J_1}\right)^2 -\frac{1}{2}\left(\S_{J_1}^{* T} E_{J_1}\right)^2
\end{align}

Therefore, the following expression can be derived as:
\begin{equation}
\begin{aligned}
\dot{V}_1  \leq&-\left(\rho_1-2\right) e_1^2-\frac{\gamma_1}{2} \tilde{\S}_{f_1}^T \tilde{\S}_{f_1}-\frac{\varepsilon_{c_1}}{2}\left(\tilde{\S}_{c_1}^T E_{J_1}\right)^2  -\frac{\varepsilon_{c_1}}{2}\left(\tilde{\S}_{a_1}^T E_{J_1}\right)^2+M_1+\frac{1}{2} e_2^2,
\end{aligned}
\end{equation}
where $M_1=\left(\frac{\varepsilon_{a_1}}{2}+\frac{\varepsilon_{c_1}}{2}\right)\left(\S_{J_1}^{* T} E_{J_1}\right)^2+\frac{\gamma_1}{2} \S_{f_1}^{* T} \S_{f_1}^*+\frac{1}{2}\dot{y}_r^2+\frac{1}{2} \omega_{f_1}^2$.

The $\chi_{\Pi_1{ }^{-1}}^{\max }$ and $\chi_{E_{J_1}}^{\min }$ are defined as the maximum and minimum eigenvalues of $\Pi_1^{-1}$ and $E_{J_1} E_{J_1}^T$, respectively, resulting in the following fact
\begin{equation}
\begin{aligned}
-\tilde{\S}_{f_1}^T \tilde{\S}_{f_1} & \leq-\frac{1}{\chi_{\Pi_1^{-1}}^{\max }} \tilde{\S}_{f_1}^T \Pi_1^{-1} \tilde{\S}_{f_1} \\
-\left(\tilde{\S}_{c_1}^T E_{J_1}\right)^2 & \leq-\chi_{E_{J_1}}^{\min } \tilde{\S}_{c_1}^T \tilde{\S}_{c_1} \\
-\left(\tilde{\S}_{a_1}^T E_{J_1}\right)^2 & \leq-\chi_{E_{J_1}}^{\min } \tilde{\S}_{a_1}^T \tilde{\S}_{a_1}
\end{aligned}
\end{equation}

From all the conclusions deduced above, we get
\begin{equation}
\begin{aligned}
\dot{V}_1 \leq&-\left(\rho_1-2\right) e_1^2-\frac{\gamma_1}{2 \chi_{\Pi_1^{-1}}^{\max }} \tilde{\S}_{f_1}^T \Pi_1^{-1} \tilde{\S}_{f_1} -\frac{\varepsilon_{c_1}}{2} \chi_{E_{J_1}}^{\min } \tilde{\S}_{c_1}^T \tilde{\S}_{c_1}-\frac{\varepsilon_{c_1}}{2} \chi_{E_{J_1}}^{\min } \tilde{\S}_{a_1}^T \tilde{\S}_{a_1} +M_1+\frac{1}{2} e_2^2.
\end{aligned}
\end{equation}

\emph{Step \textsl{i ( i = 2,...,n-1) :}} Below is the design of the Lyapunov candidate function
\begin{equation}
\begin{aligned}
V_i =&\sum_{j=1}^{i-1} V_j+\frac{1}{2} e_i^2+\frac{1}{2} \tilde{\S}_{f_i}^T\Pi_i^{-1} \tilde{\S}_{f_i}  +\frac{1}{2} \tilde{\S}_{c_i}^T \tilde{\S}_{c_i}+\frac{1}{2} \tilde{\S}_{a_i}^T\tilde{\S}_{a_i},
\end{aligned}
\end{equation}
where $V_j=\frac{1}{2} e_j^2+\frac{1}{2} \tilde{\S}_{f_j}^T \Pi_j^{-1} \tilde{\S}_{f_j}+\frac{1}{2} \tilde{\S}_{c_j}^T \tilde{\S}_{c_j}+\frac{1}{2} \tilde{\S}_{a_j}^T \tilde{\S}_{a_j} $ denotes the Lyapunov function of the ith subsystem.

Based on (26), and (35)-(37), we get
\begin{equation}
\begin{aligned}
\dot{V}_i=&\sum_{j=1}^{i-1} \dot{V}_j+e_i\left(-\rho_i e_i-\tilde{\S}_{f_i}^T E_{f_i}-\frac{1}{2} \hat{\S}_{a_i}^T E_{J_i}+e_{i+1}-\dot{\hat{\alpha}}_{i-1}^*+\omega_{f_i}\right) +\tilde{\S}_{f_i}^T\left(E_{f_i} e_i-\gamma_i \hat{\S}_{f_i}\right)
\\&-\varepsilon_{c_i} \tilde{\S}_{c_i}^T E_{J_i} E_{J_i}^T \hat{\S}_{c_i} -\tilde{\S}_{a_i}^T E_{J_i} E_{J_i}^T\left(\varepsilon_{a_i}\left(\hat{\S}_{a_i}-\hat{\S}_{c_i}\right)+\varepsilon_{c_i} \hat{\S}_{c_i}\right).
\end{aligned}
\end{equation}

The $\chi_{\Pi_i^{-1}}^{\max }$ and $\chi_{E_{J_i}}^{\min }$ are defined as the maximum and minimum eigenvalues of $\Pi_i^{-1}$ and $E_{J_i} E^T_{J_i}{ }$, respectively. Similar to steps (58)-(62), we can calculate the following expression:
\begin{equation}
\begin{aligned}
\dot{V}_i \leq&\sum_{j=1}^{i-1} \dot{V}_j-\left(\rho_i-3\right) e_i^2-\frac{\gamma_i}{2 \chi_{\Pi_i^{-1}}^{\max }} \tilde{\S}_{f_i}^T \Pi_i^{-1} \tilde{\S}_{f_i} -\frac{\varepsilon_{c_i}}{2} \chi_{E_{J_i}}^{\min } \tilde{\S}_{c_i}^T \tilde{\S}_{c_i}-\frac{\varepsilon_{c_i}}{2} \chi_{E_{J_i}}^{\min } \tilde{\S}_{a_i}^T \tilde{\S}_{a_i} +M_i+\frac{1}{2} e^2_{i+1},
\end{aligned}
\end{equation}
where $M_i=\left(\frac{\varepsilon_{a_i}}{2}+\frac{\varepsilon_{c_i}}{2}\right)\left(\S_{J_i}^{* T} E_{J_i}\right)^2+\frac{\gamma_i}{2} \S_{f_i}^{* T} \S_{f_i}^*+\frac{1}{2} {\dot{\hat{\alpha}}^{* 2}_{i-1}}+\frac{1}{2} \omega_{f_i}^2$.

\emph{Step \textsl{n :}} Considering the Lyapunov candidate function as
\begin{equation}
\begin{aligned}
V_n =&\sum_{j=1}^{n-1} V_j+\frac{1}{2} e_n^2+\frac{1}{2} \tilde{\S}_{f_n}^T \Pi_n^{-1} \tilde{\S}_{f_n} +\frac{1}{2} \tilde{\S}_{c_n}^T \tilde{\S}_{c_n}+\frac{1}{2} \tilde{\S}_{a_n}^T \tilde{\S}_{a_n}.
\end{aligned}
\end{equation}

By calculating the time derivative of $V_n$ along (38), (49)-(51), one can obtain
\begin{equation}
\begin{aligned}
\dot{V}_n=&\sum_{j=1}^{n-1} \dot{V}_j+\tilde{\S}_{f_n}^T\left(E_{f_n} e_n-\gamma_n \hat{\S}_{f_n}\right)-\varepsilon_{c_n} \tilde{\S}_{c_n}^T E_{J_n} E_{J_n}^T \hat{\S}_{c_n} -\tilde{\S}_{a_n}^T E_{J_n} E^T_{J_n}{ }\left(\varepsilon_{a_n}\left(\hat{\S}_{a_n}-\hat{\S}_{c_n}\right)+\varepsilon_{c_n} \hat{\S}_{c_n}\right)\\ 
& +e_n\left(u+f_n\left(\bar{x}_n, u\right)+\sigma\left(t-T_0\right) \lambda(x, u)-\dot{\hat{\alpha}}_{n-1}^*\right).
\end{aligned}
\end{equation}

By employing Young's inequality and Assumption 1, it yields
\begin{equation}
\begin{aligned}
e_n\left(f_n\left(\bar{x}_n, u\right)+\right. \left.\sigma\left(t-T_0\right) \lambda(x, u)\right) \leq\left|e_n\right||\bar{h}(x, u)| \leq \frac{1}{2} e_n^2 \bar{h}^2(x, u)+\frac{1}{2} \leq e_n h(x, u)+\frac{1}{2}.
\end{aligned}
\end{equation}

Clearly, where $h(x, u)=\frac{1}{2} e_n \bar{h}^2(x, u)$ denotes an unknown smooth function. Therefore, $h(x, u)$ can be approximated as
\begin{equation}
\begin{aligned}
h\left(x, u_f \mid \S_{f_n}^{* T}\right)=\S_{f_n}^{* T} E_{f_n}\left(\bar{x}_n, u_f\right)+\omega_{f_n}\left(\bar{x}_n, u_f\right),
\end{aligned}
\end{equation}
$u_f$ is defined as the filtered signal, given by
\begin{equation}
\begin{aligned}
u_f=H_L(e) u \approx u,
\end{aligned}
\end{equation}
where $H_L(e)$ represents Butterworth low-pass filer.

\emph{Remark \textsl{4:}} According to\cite{29}, it is known that using equation (16) directly in controller design can lead to algebraic loop issues. To address potential algebraic loop issues, it is a practical approach to utilize the filtered signal $u_f$ considering the low-pass characteristics commonly found in most actuators.

Based on (53) and (54), it is possible to design time-varying parameter $\left|\mu_i\right| \leq 1$, $i=1,2$ such that ${U}({t})={u}({t})+\beta\mathrm{\mu}_1 u(t)+\mathrm{\mu}_2 \theta$, leading to
\begin{equation}
\begin{aligned}
{u}({t})=\frac{{U}({t})}{1+\mu_1 \beta}-\frac{\mathrm{\mu_2} \theta}{1+\mu_1 \beta}.
\end{aligned}
\end{equation}

Since $\left|\mu_i\right| \leq 1, i=1,2$, we can conclude
\begin{equation}
\begin{aligned}
\frac{e_n {U}({t})}{1+\mu_1 \beta} \leq \frac{e_n {U}({t})}{1+\beta}, \frac{\mathrm{\mu}_2 \theta}{1+{\mu_1} \beta} \leq\left|\frac{\theta}{1-\beta}\right|.
\end{aligned}
\end{equation}

Based on (52) and (75), one has
\begin{equation}
\begin{aligned}
e_n {u}({t}) \leq&-\left|e_n \hat{u}\right|+\left|e_n\right||\hat{u}|-e_n \hat{u} \tanh \left(\frac{e_n \hat{u}}{v}\right) +\left|e_n\right||\zeta|-e_n \zeta \tanh \left(\frac{e_n \zeta}{v}\right)\\
&-\left|e_n\right||\zeta|  +\left|e_n\right|\left|\frac{\theta}{1-\beta}\right|.
\end{aligned}
\end{equation}

Defined that $\chi_{\Pi_n^{-1}}^{\max }$ and $\chi_{E_{J_n}}^{\min }$ represent the maximum and minimum eigenvalues of $\Pi_n^{-1}$ and $E_{J_n} E^T_{J_n}{ }$, respectively. According to $\left|e_n\right|-e_n \tanh \left(\frac{e_n}{v}\right) \leq \Delta v, \Delta=0.2785$ and $\zeta>\frac{\theta}{1-\beta}$, (38), (72), (73), similarly we get
\begin{equation}
\begin{aligned}
\dot{V}_n \leq& \sum_{j=1}^{n-1} \dot{V}_j-\left(\rho_n-3\right) e^2_n{ }+M_n-\frac{\varepsilon_{c_n}}{2} \chi_{E_{J_n}}^{\min} \tilde{\S}_{a_n}^T \tilde{\S}_{a_n} -\frac{\varepsilon_{c_n}}{2} \chi_{E_{J_n}}^{\min } \tilde{\S}_{c_n}^T \tilde{\S}_{c_n}  -\frac{\gamma_n}{2 \chi_{\Pi^{-1}_n}^{\max } \tilde{\S}_{f_n}^T} \Pi^{-1}_n \tilde{\S}_{f_n},\\
\end{aligned}
\end{equation}
where $ M_n=\left(\frac{\varepsilon_{a_n}}{2}+\frac{\varepsilon_{c_n}}{2}\right)\left(\S_{J_n}^{* T} E_{J_n}\right)^2+\frac{\gamma_n}{2} \S_{f_n}^{* T} \S_{f_n}^*+\frac{1}{2} {\dot{\hat{\alpha}}^{* 2}_{n-1}}+\frac{1}{2} \omega_{f_n}^2+\frac{1}{2}+$ $2 \Delta v$. Since all terms of $M_i$, $i=1,\ldots,n$ are bounded, $M_i$ can be bounded by a constant $m_i$, i.e. $M_i \leq m_i$. Let $a_i=\min \left\{2\left(\rho_i-3\right), \frac{\gamma_i}{\chi_{\Pi_i^{-1}}^{m a x}}, \varepsilon_{c_i} \chi_{E_{J_i}}^{\min }\right\}$, then we can rewrite (77) as follows
\begin{equation}
\begin{aligned}
\dot{V}_n \leq \sum_{j=1}^n\left(-a_j V_j+m_j\right).
\end{aligned}
\end{equation}

Theorem 1: Considering the nonlinear system (1), the optimal controllers (18), (34), (48), the updating laws of the identifier, critic, actor NNs (19), (35), (49), (20), (36), (50), (21), (37), and (51), along with the event-triggered controller described by (52). Under Assumption 1, if the design parameters are properly chosen, the following theoretical conclusions can be derived: all error signals are SGUUB, furthermore, the system exhibits good tracking performance, and Zeno behavior is avoided.

 Proof:
 Let ${a}=\min \left\{a_1, a_2, \ldots, a_n\right\}$ and ${m}=\sum_{i=1}^n m_i$, then (78) can become the following one:
\begin{equation}
\begin{aligned}
\dot{V}_n \leq-a V_n+m.
\end{aligned}
\end{equation}

Applying Lemma 1 to (79) yields
\begin{equation}
\begin{aligned}
V_n \leq e^{-a t} V_n(0)+\frac{m}{a}\left(1-e^{-a t}\right).
\end{aligned}
\end{equation}

Thus, one can draw the conclusion that all error signals $\tilde{\S}_{f_i}, \tilde{\S}_{a_i}, \tilde{\S}_{c_i}$ and $e_i$ are SGUUB.

In order to demonstrate the avoidance of Zeno behavior, we design the following function:
\begin{equation}
\begin{aligned}
d|{u}({t})-{U}({t})|& =\operatorname{sign}\left({U}\left({t}_k\right)-{U}({t})\right) d\left({U}\left({t}_k\right)-{U}({t})\right) \\
& \leq|d{U}({t})|.
\end{aligned}
\end{equation}

By (64), it can be inferred that $d{U}({t})$ is bounded, therefore $|d{U}({t})|$ can be bounded by a positive constant $\delta$, that is $|d{U}({t})| \leq \delta$. Then it can be derived ${u}\left({t}_k\right)-{U}\left({t}_k\right)=0$ and $\lim _{t \rightarrow {t}_{k+1}}({u}({t})-$ ${U}({t}))=\beta\left|{u}\left({t}_{k+1}\right)\right|+\theta \geq \theta$, then ${t}_{k+1}-{t}_k \geq \frac{\theta}{\delta}$. Therefore, the conclusion can be drawn to avoid Zeno behavior.

\section{Simulation}
To demonstrate the feasibility of the proposed strategy, numerical simulations are conducted as follows
\begin{equation}
\left\{
\begin{array}{lr}
\dot{x}_1 =x_{2}+f_1\left(\bar{x}_1\right) \\
\dot{x}_2 =u+f_2\left(\bar{x}_2\right)+\sigma\left(t-T_0\right) \lambda(x, u)  \\
{y}_r =sin(t)
\end{array}
\right.
\end{equation}
where $f_1\left(\bar{x}_1\right)={x}_1sin({x}_1)$ and  $f_2\left(\bar{x}_2\right)={x}_2cos({x}_1)$. Besides, the initial values are selected as ${x}_i(0)=0$, ${\S}_{f_i}(0)={\S}_{c_i}(0)={\S}_{a_i}(0)=0.3 (i=1,2)$. Additionally, the parameters are selected as: $\Pi_1=15$, $ \Pi_2=20$, $\gamma_i=3$, $\varepsilon_{c_i}=15$, $\varepsilon_{a_i}=18$, $\rho_1=40$, $\rho_2=50$, $\beta=0.2$, $\zeta=3$, $v=0.2$, $\theta=4$; Furthermore, the Butterworth low-pass filer and a time curve of a fault are selected as  $H_L(e)=1/(e^2+1.141e+1)$, $
\lambda(x, u)=4({x}_1{x}_2+sin(u))+2
$, respectively. The definition of a fault function is as follows
$$
\sigma\left(t-T_0\right)= \begin{cases}0, & t<10 \\ 1-\mathrm{e}^{-20\left(t-10\right)}, & t \geq 10\end{cases}
$$

The simulation results are shown in Figures 1-5. Fig. 1 shows the trajectories of ${x}_1$ and ${y}_r$, indicating good tracking performance. Fig. 2 displays the trajectories of the control input signal and the tracking error ${e}_1$ can converge to zero is confirmed. From Fig. 3, the boundedness of the NN weights for the identifier, critic, and actor are demonstrated. Fig. 4 displays the cost functions $c_1\left(e_1, \alpha_1\right)$ and $c_2\left(e_2, u\right)$ of two backstepping steps. Fig. 5 represents the triggering intervals $\Delta{T}={t}_{k+1}-{t}_{k}$. The control update time intervals for the event-triggered controller in this article are always larger than a time step of 0.001, effectively avoiding the occurrence of Zeno behavior. The overall triggering times are 740.

To validate the optimization performance of the proposed control method, we employed the general backstepping algorithm as described in reference \cite{37}. The comparative results are depicted in Figs. 6 and 7. As illustrated in Fig. 6, both approaches exhibit equivalent tracking performance; however, a comparison of their cost functions can be observed from Fig. 7. Through an analysis of Figs. 6 and 7, we draw a direct conclusion: under identical tracking performance conditions, the proposed control scheme demonstrates lower cost.

\begin{figure}[H]
	\centering
	\includegraphics[scale=0.6]{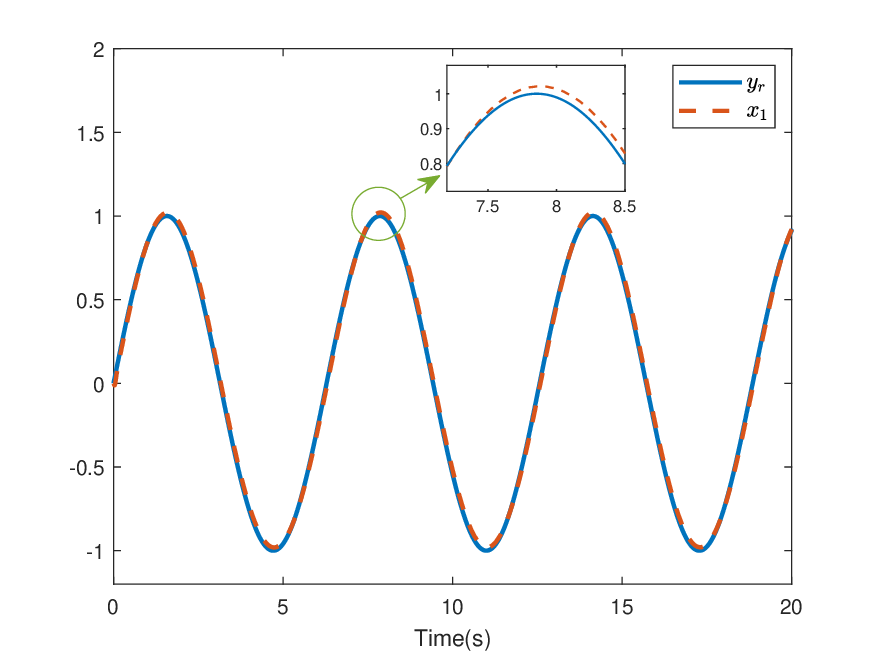}
	\caption{Trajectories of ${x}_1$ and ${y}_r$}
	\label{fig1}
\end{figure}
\begin{figure}[H]
	\centering
	\includegraphics[scale=0.6]{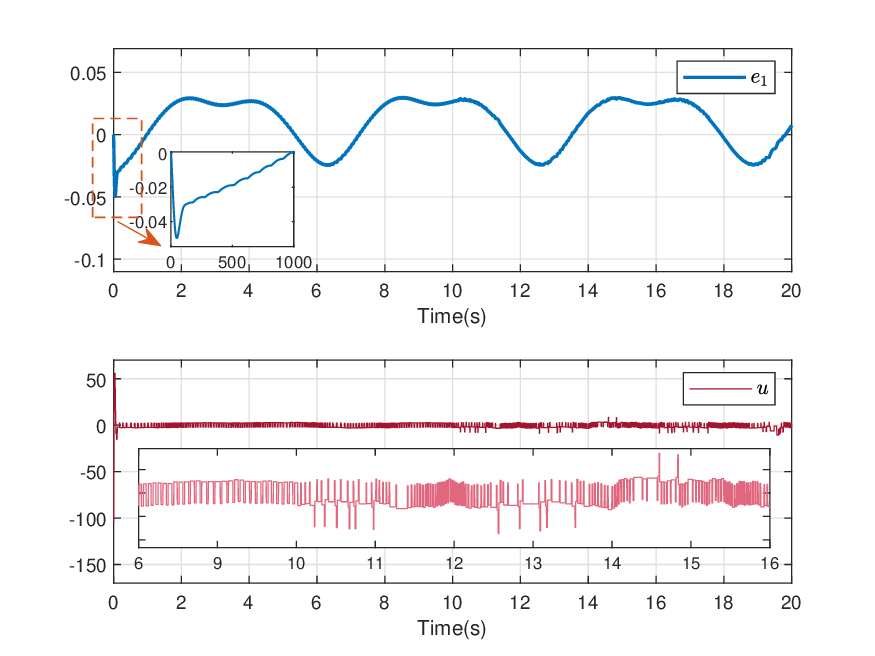}
	\caption{The tracking error $e_1$ and the actual control signal $u$}
	\label{fig2}
\end{figure}
\begin{figure}[H]
	\centering
	\includegraphics[scale=0.6]{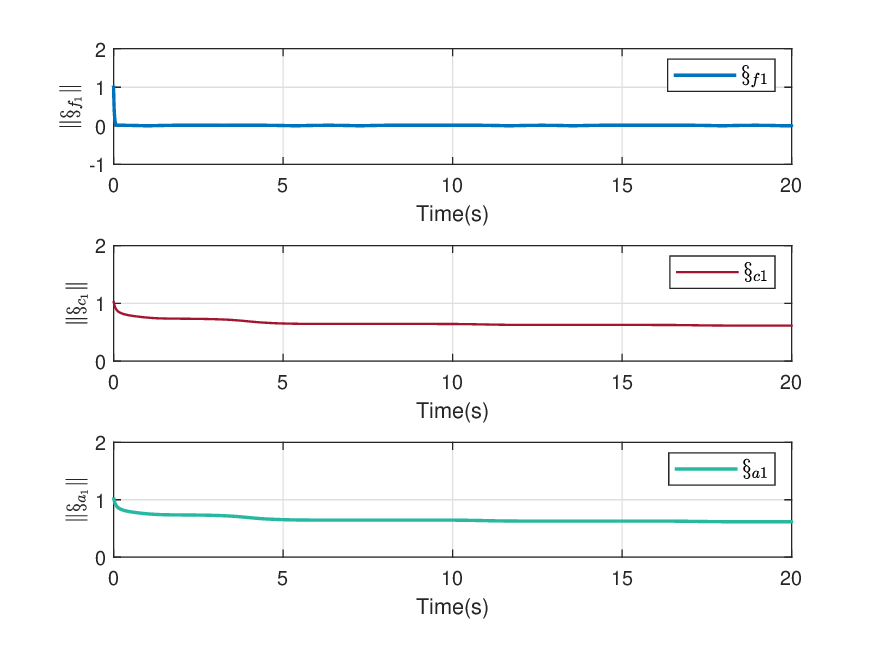}
	\caption{ $\Vert \hat{{\S}}_{f_1} \Vert$, $\Vert \hat{{\S}}_{c_1} \Vert$, $\Vert \hat{{\S}}_{a_1} \Vert$ of the first step}
	\label{fig3}
\end{figure}
\begin{figure}[H]
\vspace{-5.5cm}
	\centering
	\includegraphics[scale=0.6]{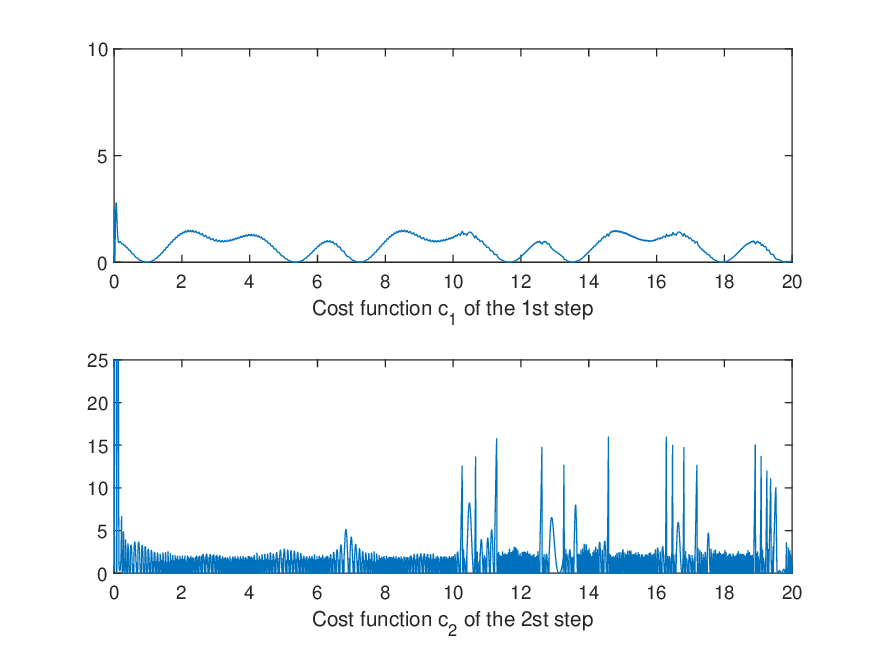}
\vspace{-5.5cm}
	\caption{The cost functions for the first and second steps}
	\label{fig4}
\end{figure}
\begin{figure}[H]
	\centering
	\includegraphics[scale=0.6]{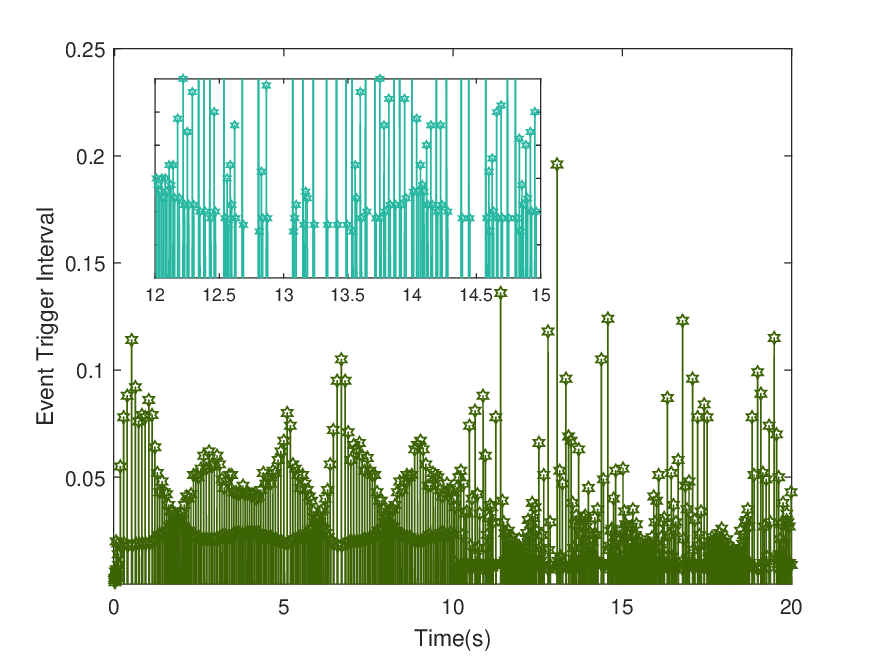}
	\caption{The triggering interval $\Delta{T}={t}_{k+1}-{t}_{k}$}
	\label{fig5}
\end{figure}
\begin{figure}[H]
\vspace{-5.5cm}
	\centering
	\includegraphics[scale=0.6]{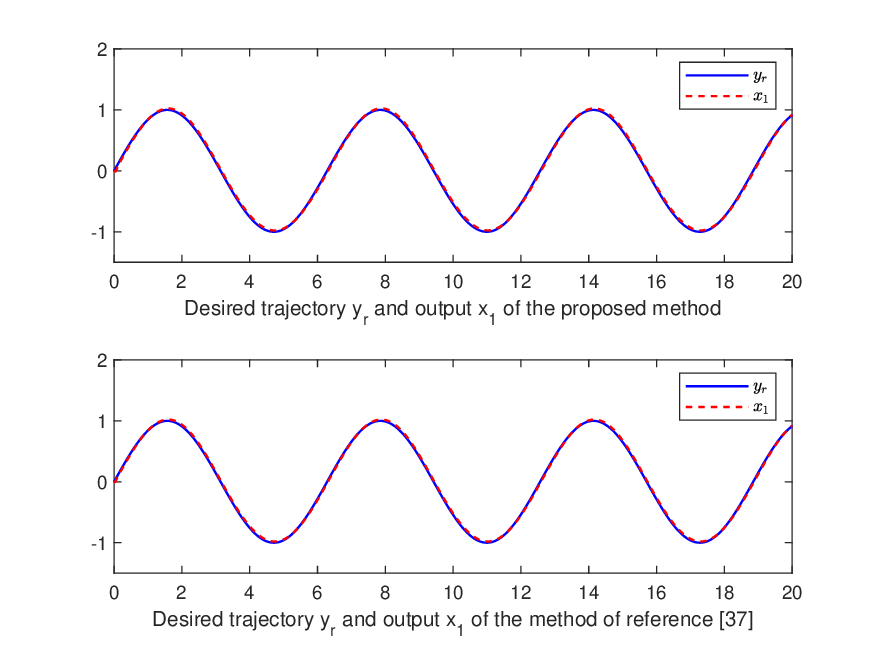}
\vspace{-5.5cm}
	\caption{Two tracking performances}
	\label{fig6}
\end{figure}
\begin{figure}[H]
\vspace{-5.5cm}
	\centering
	\includegraphics[scale=0.6]{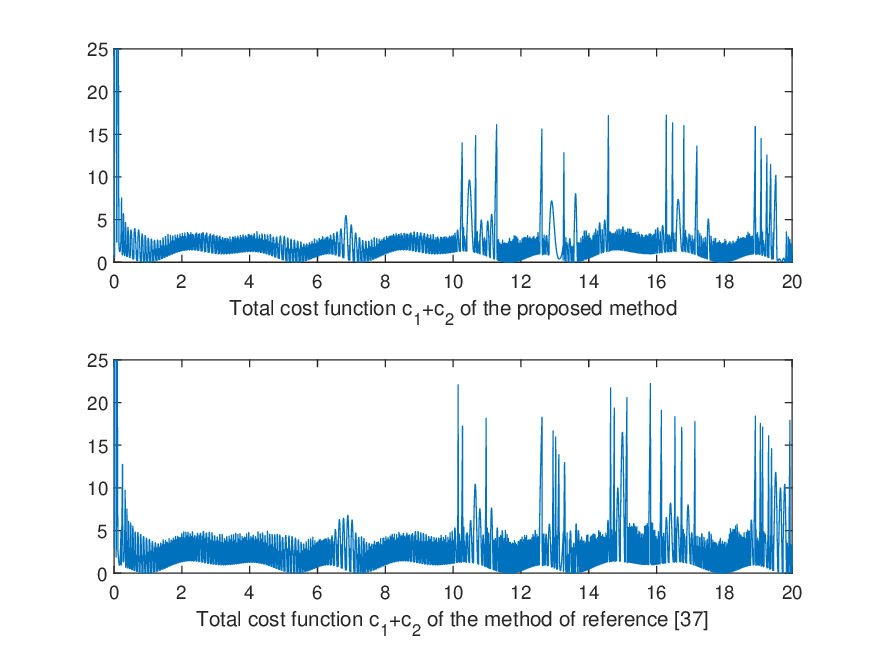}
\vspace{-5.5cm}
	\caption{Two total cost functions}
	\label{fig7}
\end{figure}

\section{Conclusion\vspace{-0.8em}}
This article proposes an event-triggered optimal tracking control scheme for a class of SFNSs with non-affine and nonlinear faults, utilizing a simplified identifier-critic-actor structure of the RL algorithm to achieve optimal control performance. The proposed simplified RL algorithm relaxes the persistent excitation condition and is designed by deriving update rules from the negative gradient of a simple positive function related to the HJBE. A FTC method is developed by applying filtered signals for controller design, and an event-triggered mechanism is adopted to reduce communication resource consumption and avoid Zeno behavior. Finally, the theoretical analysis and simulation results demonstrate the feasibility and effectiveness of the proposed scheme.\\

\noindent
\textbf{Funding:\vspace{0.4em}} This work was supported by the National and Science Foundation of China 62276214.\\
\textbf{Conflict of Interest\vspace{0.4em}} The authors declare that they have no conflict of interest.\\
\textbf{Ethical Approval\vspace{0.4em}} This article does not contain any studies with human participants or animals per-formed by any of the authors.\\
\textbf{Data availability statement}: All data generated or analyzed during this study are included in this published article.


\begin{thebibliography}{99}
        

	{\bibitem{1}
		Reshmin, S. A.
		\newblock Properties of the time-optimal control for lagrangian single-degree-of-freedom systems.
		\newblock {\em IEEE Transactions on Automatic Control}, 60(12), 3350-3355(2015).}
	
	{\bibitem{2}
		Doelman, R., Dominicus, S., Bastaits, R., and Verhaegen, M.
		\newblock Systematically structured  $H_2$  optimal control for truss-supported segmented mirrors.
		\newblock {\em IEEE Transactions on Control Systems Technology}, 27(5),  2263-2270(2019).}
	
	{\bibitem{3}
		Demirel, B., Ghadimi, E., Quevedo, D. E., and Johansson, M.
		\newblock Optimal control of linear systems with limited control actions: threshold-based event-triggered control.
		\newblock {\em IEEE Transactions on Control of Network Systems}, 5(3), 1275-1286(2018).}
	
	{{\bibitem{4}
		Van Berkel, K., Titulaer, R., Hofman, T., Vroemen, B., and Steinbuch, M.
		\newblock From optimal to real-time control of a mechanical hybrid powertrain.
		\newblock {\em IEEE Transactions on Control Systems Technology}, 23(2), 670-678(2015).}}
	
	{\bibitem{5}
		Fan, Z.-X., Li, S., and Liu, R.
		\newblock ADP-based optimal control for systems with mismatched disturbances: a PMSM application.
		\newblock {\em IEEE Transactions on Circuits and Systems II: Express Briefs}, 70(6), 2057-2061(2023).}

{\bibitem{1111}
		Liu, L., Li, Z., Chen, Y., and Wang, R.
		\newblock Disturbance observer based adaptive intelligent control of marine vessel with position and heading constraint condition related to desired output.
		\newblock {\em IEEE Transactions on Neural Networks and Learning Systems}, 1-10(2022).

	{\bibitem{6}
		Wang, D., Qiao, J., and Cheng, L.
		\newblock An approximate neuro-optimal solution of discounted guaranteed cost control design.
		\newblock {\em IEEE Transactions on Cybernetics}, 52(1), 77-86(2022).}
	
	
	{\bibitem{7}
		Werbos, P. J.
		\newblock Approximate dynamic programming for real-time control and neural modeling.
		\newblock {\em Handbook Intell. Control Neural Fuzzy Adaptive Approaches}, 15, 493-525(1992).}
	
	
	{\bibitem{8}
		Lin, Ziyu., Ma, Jun., Duan, Jingliang., Li, Shengbo Eben., Ma, Haitong., Cheng, Bo., and Lee, Tong Heng.
		\newblock Policy iteration based approximate dynamic programming toward autonomous driving in constrained dynamic environment.
		\newblock {\em IEEE Transactions on Intelligent Transportation Systems}, 24(5), 5003-5013(2023).}}
	
	{\bibitem{9}
		Shi, L., Wang, X., and Cheng, Y.
		\newblock Afe reinforcement learning-based robust approximate optimal control for hypersonic flight vehicles.
		\newblock {\em IEEE Transactions on Vehicular Technology}, 72(9), 11401-11414(2023).}
	
	{\bibitem{10}
		An, T., Wang, Y., Liu, G., Li, Y., and Dong, B.
		\newblock Cooperative game-based approximate optimal control of modular robot manipulators for human--robot collaboration.
		\newblock {\em IEEE Transactions on Cybernetics}, 53(7), 4691-4703(2023).}
	
	{\bibitem{11}
		Pecioski, D., Gavriloski, V., Domazetovska, S., and Ignjatovska, A.
		\newblock  An overview of reinforcement learning techniques.
		\newblock {\em 2023 12th Mediterranean Conference on Embedded Computing (MECO)}, 1-4(2023).}

{\bibitem{2222}
                  Liu, L., Cui, Y., Liu, Y.-J., and Tong, S.
                  \newblock Observer based adaptive neural output feedback constraint controller design for switched systems under average dwell time.
                 \newblock {\em IEEE Transactions on Circuits and System I: Regular Papers}, 68(9), 3901-3912(2021).}
	
	{\bibitem{13}
		Liu, Y., Zhu, Q., and Wen, G.
		\newblock Adaptive tracking control for perturbed strict-feedback nonlinear systems based on optimized backstepping technique.
		\newblock {\em IEEE Transactions on Neural Networks and Learning Systems}, 33(2), 853-865(2022).}

	{\bibitem{14}
		Wang, H., and Bai, W.
		\newblock Finite-time adaptive fault-tolerant control for strict-feedback nonlinear systems.
		\newblock {\em 2019 Chinese Control And Decision Conference (CCDC)}, Nanchang, China, 5200-5204(2019).}

	{\bibitem{15}
		Pang N., Wang X., and Wang Z. M.
		\newblock Event-triggered adaptive control of nonlinear systems with dynamic uncertainties: The switching threshold case.
		\newblock {\em IEEE Transactions on Circuits and Systems II: Express Briefs}, 69(8), 3540-3544(2022).}

{\bibitem{3333}
		Wang, F., Liu, Z., Zhang, Y., and Chen, C. L. P.
		\newblock Adaptive fuzzy control for a class of stochastic pure-feedback nonlinear systems with unknown hysteresis.
		\newblock {\em IEEE Transactions on Fuzzy Systems}, 24(1), 140-152(2016).}

{\bibitem{4444}
		Bhasin, S., Kamalapurkar, R., Johnson, M., Vamvoudakis, K. G., Lewis, F. L., and Dixon, W. E.
		\newblock A novel actor–critic-identifier architecture for approximate optimal control of uncertain nonlinear systems.
		\newblock {\em Automatica}, 49(1), 82-92(2013).}	

	{\bibitem{16}
		Wen, G., Ge, S. S., and Tu, F.
    \newblock Optimized backstepping for tracking control of strict-feedback systems.
    \newblock {\em IEEE Transactions on Neural Networks and Learning Systems}, 29(8), 3850-3862(2018).}
	
	{\bibitem{17}
		Wen, G., Xu, L., and Li, B.
		\newblock Optimized backstepping tracking control using reinforcement learning for a class of stochastic nonlinear strict-feedback systems.
		\newblock {\em  IEEE Transactions on Neural Networks and Learning Systems}, 34(3), 1291-1303(2023).}
	
	{\bibitem{18}
		Wang, X., Guang, W., Huang, T., and Kurths, J.
		\newblock Optimized adaptive finite-time consensus control for stochastic nonlinear multiagent systems with non-affine nonlinear faults.
		\newblock {\em  IEEE Transactions on Automation Science and Engineering}, 1-12(2023).}
	
         {\bibitem{5555}
    Pang, N., Wang, X., and Wang, Z. M.
    \newblock Observer-based event-triggered adaptive control for nonlinear multiagent systems with unknown states and disturbances
    \newblock {\em IEEE Transactions on Neural Networks and Learning Systems}, 34(9), 6663-6669(2021).}


         {\bibitem{6666}
    Wang, X., Xu, R., Huang, T. W., Kurths, J.
    \newblock Event-triggered adaptive containment control for heterogeneous stochastic nonlinear multiagent systems.
    \newblock {\em IEEE Transactions on Neural Networks and Learning Systems}, 1-11(2023).}

	{\bibitem{19}
		Dong, Y., and Lin, Z.
		\newblock An event-triggered observer and its applications in cooperative control of multiagent systems.
		\newblock {\em  IEEE Transactions on Automatic Control}, 67(7), 3647-3654(2022).}
	
	{\bibitem{20}
		Yu, H., Hao, F., and Chen, T.
		\newblock A uniform analysis on input-to-state stability of decentralized event-triggered control systems.
		\newblock {\em  IEEE Transactions on Automatic Control}, 64(8), 3423-3430(2019).}
	
	{\bibitem{21}
		Delimpaltadakis, Georgios., and Mazo, Manuel.
		\newblock Abstracting the traffic of nonlinear event-triggered control systems.
		\newblock {\em  IEEE Transactions on Automatic Control}, 68(6), 3744-3751(2023).}
	
	
	{\bibitem{22}
		Wang, Z. C., and Wang, X.
		\newblock Event-Triggered Containment Control for Nonlinear Multiagent Systems via Reinforcement Learning.
		\newblock {\em  IEEE Transactions on Circuits and Systems II: Express Briefs}, 70(8), 2904-2908(2023).}

        
       {\bibitem{7777}
    Wang, Z., Wang, H., Wang, X., Pang, N., and Shi, Q.
    \newblock Event-triggered adaptive neural control for full state-constrained nonlinear systems with unknown disturbances.
    \newblock {\em Cognitive Computation}, 16, 717-726(2024).}

         {\bibitem{8888}
    Wang, Z., Wang, X., and Pang, N.
    \newblock Adaptive fixed-time control for full state-constrained nonlinear systems: switched-self-triggered case.
    \newblock {\em IEEE Transactions on Circuits and Systems II: Express Briefs}, 71(2), 752-756(2024).}


	
	{\bibitem{23}
		Ofodile, Ikechukwu., Ofodile-Keku, Nkemdilim., Jemitola, Paul., Anbarjafari, Gholamreza., and Slavinskis, Andris.
		\newblock Integrated anti-windup fault-tolerant control architecture for optimized satellite attitude stabilization.
		\newblock {\em  IEEE Journal on Miniaturization for Air and Space Systems}, 2(4), 189-198(2021).}
	
	{\bibitem{24}
		Shen, Q., Wang, D., Zhu, S., and Poh, E. K.
		\newblock Integral-type sliding mode fault-tolerant control for attitude stabilization of spacecraft.
		\newblock {\em  IEEE Transactions on Control Systems Technology}, 23(3), 1131-1138(2015).}
	
	{\bibitem{25}
		Tang, H., Chen, Y., and Zhou, A.
		\newblock Actuator fault-tolerant control for four-wheel-drive-by-wire electric vehicle.
		\newblock {\em  IEEE Transactions on Transportation Electrification}, 8(2), 2361-2373(2022).}
	
	{\bibitem{26}
		Wang, Z. M.
		\newblock Hybrid Event-triggered Control of Nonlinear System with Full State Constraints and Disturbance.
		\newblock {\em  2024, arXiv:2405.13564.} [Online]. Available: https://arxiv.org/abs/2405.13564}
	
	{\bibitem{27}
		Li, X. -J., Shi, C. -X., and Yang, G. -H.
		\newblock Observer-based adaptive output-feedback fault-tolerant control of a class of complex dynamical networks.
		\newblock {\em  IEEE Transactions on Systems, Man, and Cybernetics: Systems}, 48(12), 2407-2418(2018).}
	
	{\bibitem{28}
		Yuwei, C., Aijun, L., and Xianfeng, M. 
		\newblock A fault-tolerant control method for distributed flight control system facing wing damage.
		\newblock {\em  Journal of Systems Engineering and Electronics}, 32(5), 1041-1052(2021).}
	
	{\bibitem{29}
		Sun, K., Liu, L., Qiu, J., and Feng, G.
		\newblock Fuzzy adaptive finite-time fault-tolerant control for strict-feedback nonlinear systems.
		\newblock {\em  IEEE Transactions on Fuzzy Systems}, 29(4), 786-796(2021).}
	
	{\bibitem{37}
		Ge, S. S., and Wang, C.
		\newblock Direct adaptive NN control of a class of nonlinear systems.
		\newblock {\em  IEEE Transactions on Neural Networks}, 13(1), 214-221(2002).}
	
	
	
	
	
	
	
	
	
	
	
	
	
	
	
\end{thebibliography}
\end{document}